\documentclass{svproc}
\pdfoutput=1 % force PDFLaTeX
\usepackage{hyperref,breakurl}

\usepackage{eucal}
\usepackage{stackrel}
\usepackage{mathrsfs}
\usepackage{amsmath,amssymb}
\usepackage{graphicx}
\newcommand{\Com}{\mathbb{C}}
\newcommand{\Real}{\mathbb{R}}

\newtheorem{Theorem}{Theorem}
\newcommand{\E}{\ensuremath{\mathrm{e}}}
\newcommand{\cf}{\emph{cf.}}
\newcommand{\ie}{\emph{i.e.}}
\newcommand{\eg}{\emph{e.g.}}
\newcommand{\Bounded}{\mathscr{B}}
\newcommand{\sgn}{\mathrm{sgn}}
\newcommand{\curl}{\mathop{\mathrm{curl}}\nolimits}
\newcommand{\dist}{\mathop{\mathrm{dist}}\nolimits}
\begin{document}
\mainmatter               
\title{The virial theorem and the method of multipliers in spectral theory}
\titlerunning{The virial theorem and the method of multipliers}   
\author{Lucrezia Cossetti\inst{1} \and David Krej\v{c}i\v{r}\'ik\inst{2}}  
\authorrunning{L.~Cossetti and D.~Krej\v ci\v r\'ik}  
\institute{Ikerbasque \& Universidad del Pa\'is Vasco/Euskal 
Herriko Unibertsitatea, UPV/EHU,
Aptdo.~644, 48080, Bilbao, Spain; 
\email{lucrezia.cossetti@ehu.eus} 
\smallskip \\
\and
Department of Mathematics, Faculty of Nuclear Sciences and 
Physical Engineering, Czech Technical University in Prague, 
Trojanova 13, 12000 Prague 2, Czechia;
\email{david.krejcirik@fjfi.cvut.cz}
\medskip \\
{\small 16 July	 2024}}
\maketitle             
\begin{abstract}
We provide a link between the virial theorem in functional analysis
and the method of multipliers in theory of partial differential equations.
After giving a physical insight into the techniques,
we show how to use them to deduce the absence of eigenvalues
and other spectral properties of electromagnetic quantum Hamiltonians.
We focus on our recent developments in non-self-adjoint settings,
namely on Schr\"o\-ding\-er operators with matrix-valued potentials,
relativistic operators of Pauli and Dirac types,
and complex Robin boundary conditions.
\keywords{virial theorem, method of multipliers, absence of eigenvalues,
uniform resolvent estimates, electromagnetic Schr\"odinger
and Dirac operators, non-self-adjoint perturbations,
Robin boundary conditions}
\end{abstract}
\section{Introduction}
As the final conference~\cite{Braga} amply demonstrates, 
our action has been a successful cooperation 
of various research groups  in diverse areas of mathematics. 
In this contribution, we shall reveal a close connection
between two apparently distinct tools in functional analysis 
and partial differential equations,
namely the virial theorem and the method of multipliers, respectively.
Our focus will be on surveying new developments
which has led to important applications in spectral theory
of non-self-adjoint operators during the years of the action. 

As a prelude, in Section~\ref{Sec.classical},
we start with the virial theorem in classical physics.
The formalism of quantum mechanics is recalled in Section~\ref{Sec.quantum}.
These physical concepts are not necessary for the mathematics 
considered in this paper, but they provide a useful insight,
intuition and motivation for the rigorous theorems. 
The virial theorem and the method of multipliers
are presented in Sections~\ref{Sec.virial} 
and~\ref{Sec.multiplier}, respectively.
Section~\ref{Sec.free} demonstrates the two approaches 
on the case of the free Hamiltonian,  
while electromagnetic perturbations are considered in Section~\ref{Sec.electro}.
In Section~\ref{Sec.NSA}, we present the necessary developments
to include possibly non-self-adjoint perturbations.
In Section~\ref{Sec.uniform}, it is shown how these elaborations
can be used to establish uniform resolvent estimates. 
Sections~\ref{Sec.Dirac} and~\ref{Sec.boundary} are devoted 
to applications to relativistic operators in quantum mechanics
and to boundary perturbations of the Laplacian, respectively.

This presentation is based on invited talks of the second author
in Paris~\cite{Paris}, Brijuni~\cite{Brijuni} and Bilbao~\cite{Bilbao}. 
He is grateful to the organisers for these most stimulating conferences.
Finally, as a member of the management committee 
for the Czech Republic, he should like to thank 
Marjeta Kramar Fijav\v{z} and Ivica Naki\'c
for their immense work with the action.

\section{Classical physics}\label{Sec.classical}
As a warm-up, let us recall the virial theorem in classical physics.
Consider one particle of mass~$m$ moving in the Euclidean
space~$\Real^d$ of dimension $d \geq 1$, 
subject to a force $F = -\nabla V$ 
with smooth potential $V:\Real^d \to \Real$.
The time evolution of the particle position~$x$
is given by Newton's law 
\begin{equation}\label{Newton}
  p' = F 
  \,, \qquad \mbox{where} \qquad 
  p := m x'
\end{equation}
is the momentum of the particle
and the dash denotes the time derivative. 
The total energy of the system is represented 
by the Hamiltonian function 
\begin{equation}\label{energy}
  H_V := H_0 + V
  \,, \qquad \mbox{where} \qquad 
  H_0 := \frac{p^2}{2m}
\end{equation}
is the kinetic energy (\emph{vis viva}) of the particle.
Here $p^2 := p \cdot p$, where the dot denotes the scalar product in~$\Real^d$

Now, let us ask the following question:
\begin{center}
Which potentials $V$ do not bound the particle?
\end{center}
By this we mean that the \emph{particle propagates},
which we mathematically interpret as the following divergence
of the radial velocity
\begin{equation}\label{divergence}
  T_0(t) := x(t) \cdot p(t) \xrightarrow[t \to \infty]{} \infty
  \,.
\end{equation}
Indeed, note that~\eqref{divergence} implies 
$x(t)^2 \to \infty$ as $t \to \infty$.

Following what we have learnt at school
(``if you do not know what to do, take a derivative''),
we differentiate the quantity~$T_0$ using~\eqref{Newton}
and find out the \emph{virial identity} (the first equality)
\begin{equation}\label{virial}
  T_0' =  2 H_0 - x \cdot \nabla V = \{T_0,H_V\} 
  \,.
\end{equation}
Here the Poisson bracket is defined as usual by
$$
  \{T_0,H_V\} := 
  \sum_{j=1}^d
  \left(
  \frac{\partial T_0}{\partial{x_j}}  \frac{\partial{H_V}}{\partial{p_j}}
  - \frac{\partial T_0}{\partial{p_j}}  \frac{\partial{H_V}}{\partial{x_j}}
  \right)
  ,
$$
but at this moment it is enough to consider $\{T_0,H_V\}$ 
as a shortcut for the right-hand side of~\eqref{virial}.
The following observation is an immediate consequence of 
the virial identity~\eqref{virial}.
\begin{Theorem}[the virial theorem in classical physics]\label{Thm.classical}
If $\{T_0,H_V\} \geq a$ with a positive constant~$a$
(independent of time),
then the particle propagates. 
\end{Theorem}

\begin{figure}[h!]
\begin{center}
\includegraphics[width=0.8\textwidth]{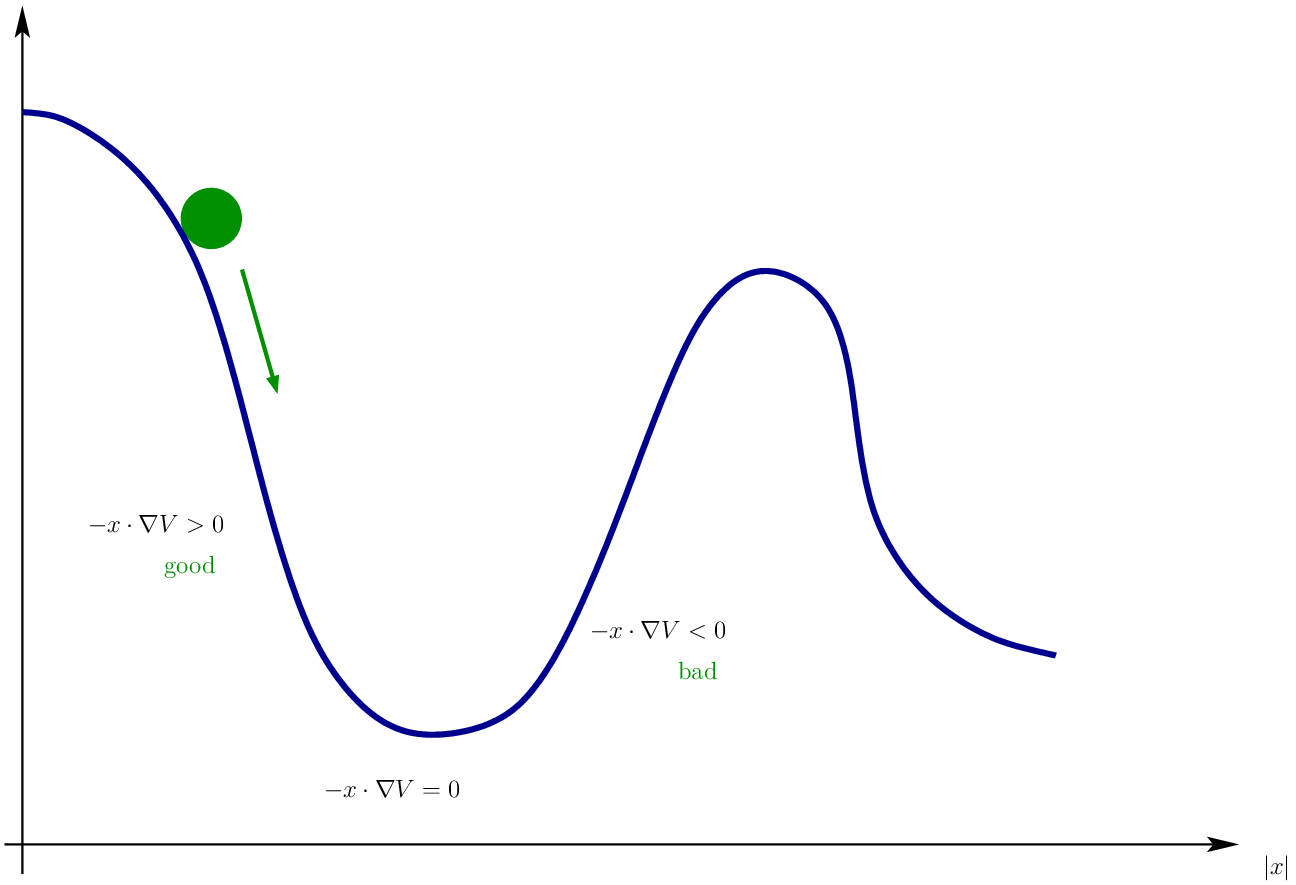}
\end{center}
\caption{A visualisation of 
the repulsivity condition~\eqref{repulsive.classical}.}
\label{Fig.chain}
\end{figure}

Since the kinetic energy~$H_0$ is a non-negative function, 
a sufficient condition to ensure the positivity $\{T_0,H_V\} \geq a$
is the repulsivity condition
\begin{equation}\label{repulsive.classical}
  - x \cdot \nabla V \geq a 
  \,.
\end{equation}
This is very intuitive if you think of the particle as
a ball moving on a mountain chain, see Figure~\ref{Fig.chain}.
More generally, the radial derivative $- x \cdot \nabla V$ 
can be zero (mountain valley) or negative (ascend)
at the times when it has a sufficiently 
large kinetic energy~$H_0$ (to escape the valley
and surmount the next peak).

\section{Quantum mechanics}\label{Sec.quantum}
In quantum mechanics, 
physical states and observables are 
represented by vectors and self-adjoint operators
in a complex (separable) Hilbert space~$\mathcal{H}$, respectively.
The expectation value of an observable~$T$ 
to be in a state~$\Psi$
is given by the inner product $\langle T \rangle := (\Psi, T\Psi)$
and the outcomes of measuring are the spectrum of~$T$.
The most prominent observable is the Hamiltonian~$H$
representing the total energy of the system.
It determines the time evolution of states  
through the Schr\"odinger equation
\begin{equation}\label{Schrodinger}
  i \, \Psi' = H\Psi \,.        
\end{equation}

Differentiating the expectation value of~$T$ with respect to time
and using~\eqref{Schrodinger},
we (at least formally) get
\begin{align}\label{A-evolution}
  \langle T \rangle'
  = \big\langle i[H,T] \big\rangle
  \,,
\end{align}
where $[H,T] := HT-TH$.
Hence the evolution of the expectation value of~$T$
is given by the expectation value of the commutator with~$H$
multiplied by~$i$ (without this multiplication,
the commutator $[H,T]$ is actually skew-adjoint).
The identity~\eqref{A-evolution} is an analogue of 
the classical formula (\cf~\eqref{virial}),
where the derivative of observables is determined 
by its Poisson bracket with the Hamiltonian.

The set of eigenvalues $\lambda \in \sigma_\mathrm{p}(H)$ 
of the Hamiltonian~$H$ are energies of the system
for which~\eqref{Schrodinger} admits stationary solutions
of the type $\E^{-i \lambda t} \psi$,
where~$\psi$ is an eigenvector of~$H$ corresponding to~$\lambda$;
it is customarily called a \emph{bound state} 
(or \emph{trapped mode}).
On the other hand, 
the continuous spectrum $\sigma_\mathrm{c}(H)$ of~$H$  
can be interpreted as energies corresponding 
to \emph{scattering states} (or \emph{propagation modes}). 
The exclusion of eigenvalues 
(specifically those embedded in the essential spectrum of~$H$)
constitutes a first step in justifying transport for a quantum system.

\section{The virial theorem}\label{Sec.virial}
How to achieve the absence of eigenvalues of a given operator~$H$?
A powerful tool is represented by an abstract version
of the \emph{virial theorem}.
Let us present a formal statement first.
We are inspired by Theorem~\ref{Thm.classical} and the fact that 
the Poisson brackets are replaced by the commutators in quantum theory. 

In addition to the self-adjoint operator~$H$,
let us consider another self-adjoint operator~$T$ in~$\mathcal{H}$.
Assume that the commutator of~$T$ with~$H$ is \emph{positive} in a sense.
For instance, in a very restrictive sense, 
that there exists a positive number~$a$ such that
(we do not care about operator domains for a moment)
\begin{equation}\label{Mourre.formal}
  i [H,T] \geq a \, I
\end{equation}
in the sense of quadratic forms in~$\mathcal{H}$.

Now, let~$\lambda$ be an eigenvalue of~$H$
corresponding to an eigenvector~$\psi$,
normalised to~$1$ in~$\mathcal{H}$. 
That is, the stationary Schr\"odinger equation
\begin{equation}\label{ev}
  H\psi = \lambda \psi
\end{equation}
holds.
Then we get a contradiction
\begin{equation}\label{contradiction}
\begin{aligned}
  a &\leq (\psi,i[H,T]\psi)
  \\
  &= i (H\psi,T\psi) - i (T\psi,H\psi)
  \\
  &= i (\lambda\psi,T\psi) - i (T\psi,\lambda\psi)
  \\
  &= 0 \,,
\end{aligned}  
\end{equation}
where the first and last equalities employ 
the self-adjointness of~$H$ and~$T$.  
Note that our convention is that the inner product $(\cdot, \cdot)$
of~$\mathcal{H}$ is linear in the second component.
Hence, the positivity of the commutator~\eqref{Mourre.formal} 
prevents the existence of eigenvalues.

To make the argument above rigorous, 
it is important to ensure that the eigenvectors of~$H$
belong to the domain of~$T$.
More universally, the conclusion can be substantially generalised
by taking a suitable definition of the commutator.
To this purpose, it is customarily assumed that~$H$ 
is of class of~$C^1(T)$, meaning that the map
$
  t \mapsto \E^{it T} (H-z)^{-1} \E^{-it T} 
$
is of class~$C^1$ for the strong topology of $\Bounded(\mathcal{H})$
for some $z \in \rho(H)$.
Then the commutator is defined via the derivative of the map at $t=0$
and the following abstract theorem holds.

\begin{Theorem}[the abstract virial theorem]\label{Thm.virial}
Let $T$ and $H \in C^1(T)$ be self-adjoint operators.
If~$\psi$ is an eigenvector of~$H$, 
then $(\psi,i[H,T]\psi)=0$.
In particular,
if~\eqref{Mourre.formal} holds with a positive constant~$a$,
then $\sigma_\mathrm{p}(H) = \varnothing$. 
\end{Theorem}

We refer to~\cite[Prop.~7.2.10]{ABG} 
for a proof of this abstract version of the virial theorem.  
For a history of the virial theorem,
see \cite[Sec.~13 \& Notes]{RS4}.

\section{The method of multipliers}\label{Sec.multiplier}
The moral of this paper is that
the virial theorem is closely related to
the \emph{method of multipliers},
usually attributed to the original 
development of Morawetz~\cite{Morawetz_1968}.

In its simplest version, the method of multipliers reads as follows.
Take an inner product of both sides of~\eqref{ev}
with the vector $\phi := i T\psi$
(this is the multiplier of the method)
and take twice the real part of the obtained identity:
\begin{equation*} 
\begin{aligned}
  (\psi,i[H,T]\psi) 
  &= (iT\psi,H\psi) + (H\psi,iT\psi) 
  \\
  &= 2\Re (\phi,H\psi) 
  \\
  &\stackrel[]{\downarrow}{=} \lambda \, 2\Re (\phi,\psi)
  \\
  &= \lambda \, [(iT\psi,\psi) + (\psi,iT\psi)] 
  \\
  &= 0
\end{aligned} 
\end{equation*}
(here the arrow points to the initial identity,
the other equalities are manipulations).
In this way we have arrived at the same identity as in~\eqref{contradiction}
and the same contradiction under the positivity hypothesis~\eqref{Mourre.formal}.

Apart from certain mathematical justifications of the manipulations above,
we see that the method of multipliers is just equivalent to the virial theorem.
However, the former turns out to be more flexible.
In particular, if~$H$ is allowed to be non-self-adjoint,    
as we shall see later.

The weak point of the hypothesis~\eqref{Mourre.formal}
is that its consequences are too strong. 
In applications, it is typically not needed to prove 
the total absence of eigenvalues.
Weaken commutator estimates, 
localised in the spectrum and with a relaxed positivity,
are usually associated with the name of Mourre~\cite{Mourre};
we refer to the book~\cite{ABG} for a comprehensive account of the theory.
Contrary to what one can occasionally read in some papers,
the method of multipliers is equally adaptable to disprove the existence 
of eigenvalues in separate subregions of the complex plane as well.

\section{The free Hamiltonian}\label{Sec.free}
Mathematically, we understand why 
the positivity of the commutator~\eqref{Mourre.formal}
is related to the (total) absence of eigenvalues of~$H$.
But, how to choose the auxiliary (so-called \emph{conjugate}) operator~$T$?
To answer this pertinent question, 
it is useful to get a physical insight first.

To this purpose, let us focus on the Hamiltonian 
of a free (\ie\ no forces) 
non-relativistic (\ie\ no spin) particle. 
It is customarily represented by the operator
\begin{equation}\label{free}
  H_0 := -\Delta
  \qquad \mbox{in} \qquad
  L^2(\Real^d) \,, 
\end{equation}
which is self-adjoint provided its domain 
is chosen to be the Sobolev space $W^{2,2}(\Real^d)$.
Note that $H_0 = p^2 = p \cdot p$, 
where $p := -i\nabla$, with domain 
being the Sobolev space~$W^{1,2}(\Real^d)$,
represents the momentum of the particle.
Accepting this definition of~$p$, 
we see that~$H_0$ coincides with
the classical kinetic energy~\eqref{energy} 
with the choice $m = \frac{1}{2}$
(without loss of generality). 
In this representation, the position of the particle
is represented by the maximal operator of multiplication 
by the space variable~$x$.

Now, let~$T_0$ be the quantum counterpart of 
the radial momentum~\eqref{divergence}:
\begin{equation}\label{dilation}
  T_0
  := \frac{x \cdot p + p \cdot x}{2}
  = -i \, x \cdot \nabla - i \, \frac{d}{2}
  \,.
\end{equation}
Note that we had to take a symmetrised version of~\eqref{divergence}
(in order to make~$T_0$ self-adjoint, at least formally),
since the observables~$x$ and~$p$ do not commute in quantum mechanics.
Then the positivity~\eqref{Mourre.formal} 
with the help of~\eqref{A-evolution} implies 
that the expectation value $\langle T_0 \rangle$ diverges for large times,
in analogy with the classical requirement~\eqref{divergence}.
It can be interpreted in physical terms
as that the particle escapes to infinity of~$\Real^d$ for large times
(for the radial derivative diverges).
That is, the particle is not bound, it propagates.
More specifically, the stationary solutions of 
the Schr\"odinger equation~\eqref{Schrodinger},
corresponding to initial data being eigenfunctions, do not exist.
These heuristic considerations suggest that~$T_0$
should be the right choice for the conjugate operator. 

It remains to analyse the validity of~\eqref{Mourre.formal} 
for the free Hamiltonian~\eqref{free}
and the radial momentum~\eqref{dilation}.
It is easily verified that (still formally) 
\begin{equation}\label{commutator.free}
  i[H_0,T_0] = 2H_0 \,.
\end{equation}
Here the right-hand side is non-negative because,
by an integration by parts,
\begin{equation}\label{positive}
  (\phi,H_0\phi) = (\phi,-\Delta\phi) = \|\nabla\phi\|^2 \geq 0
\end{equation}
for every $\phi \in W^{2,2}(\Real^d)$.
However, it is not positive in the strict sense~\eqref{Mourre.formal}
for $\sigma(H_0) = [0,\infty)$.
Nonetheless, a contradiction in the spirit of~\eqref{contradiction}
is still in order:
\begin{equation}\label{contradiction.bis}
\begin{aligned}
  2 \, \|\nabla\psi\|^2 
  &= (\psi,2H_0 \psi) 
  \\
  &\stackrel[]{\downarrow}{=} (\psi,i[H_0,T_0]\psi) 
  \\
  &= 0 
  \,,
\end{aligned} 
\end{equation}
whenever~$\psi$ is an eigenfunction of~$H_0$.
Indeed, from this identity we deduce that~$\psi$ is constant,
which is not possible for a non-trivial function in $L^2(\Real^d)$.
Therefore, we conclude with the following result.

\begin{Theorem}\label{Theorem.free}
$\sigma_\mathrm{p}(H_0) = \varnothing$.
\end{Theorem}

Indeed, this is precisely what we have just ``proved''.
In quotation marks, because there are certainly
a number of formal manipulations in the arguments given above.
First of all, the domain of the conjugate operator~$T_0$, 
formally introduced in~\eqref{dilation}, should be specified.
The customary approach is to introduce~$T_0$
as the infinitesimal generator of the dilation group~$W_t$
defined by $(W_t\psi)(x) := \E^{td/2} \psi(\E^t x)$. 
Proving that $H_0 \in C^1(T_0)$, 
Theorem~\ref{Thm.virial} and~\eqref{commutator.free}
then indeed imply Theorem~\ref{Theorem.free}.
Instead of following this direction,
we give a proof of Theorem~\ref{Theorem.free}
by means of the method of multipliers.

\begin{proof}[Theorem~\ref{Theorem.free} via the method of multipliers]
The eigenvalue equation~\eqref{ev} for the free Ha\-miltonian
precisely means that there exists a non-trivial function 
$\psi \in W^{2,2}(\Real^d)$ such that 
\begin{equation}\label{weak.formulation}
  \forall \phi \in W^{1,2}(\Real^d)
  \,, \qquad
  (\nabla\phi,\nabla\psi) = \lambda \, (\phi,\psi)
  \,.
\end{equation}
This is just a weak formulation of 
the stationary Schr\"odinger equation in~$\Real^d$.

First of all, notice that we may restrict to $\lambda \geq 0$
due to the self-adjointness of~$H_0$ and~\eqref{positive}.
In other words, the existence of non-real and negative eigenvalues
is easily disproved.

Following the formal arguments given above,
our aim is to choose~$iT_0\psi$ 
for the test function (the multiplier)~$\phi$,
where the conjugate operator~$T_0$ is given by~\eqref{dilation}. 
However, it is not clear that $\psi$ belongs to the domain of~$T_0$
(the domain of~$T_0$ has not been even discussed)
and, even if so, that $\phi \in W^{1,2}(\Real^d)$.
Indeed, the problem is the unbounded position operator~$x$ 
in the definition of~$T_0$.

To proceed rigorously, we therefore choose the regularised multiplier
\begin{equation}\label{regularised}
  \phi := x \cdot \nabla (\xi_n\psi) + \frac{d}{2} \psi
  \,,
\end{equation}
where~$\xi_n$ is the cut-off function satisfying, for every $n > 0$, 
$\xi_n(x) := \xi(x/n)$, where $\xi \in C_0^\infty(\Real^d)$
is such that $0 \leq \xi \leq 1$, $\xi(x)=1$ for every $|x| \leq 1$  
and $\xi(x)=0$ for every $|x| \geq 2$.  
Then $\phi \in W^{1,2}(\Real^d)$ because $\psi \in W^{2,2}(\Real^d)$
and the multiplication by~$x$ is bounded on the support of~$\xi_n\psi$. 
Then we get the ultimate identity $\|\nabla\psi\| = 0$
of~\eqref{contradiction.bis} after taking the limit $n \to \infty$.
\hfill\qed
\end{proof}

The specialty of the free Hamiltonian~$H_0$ is that 
we \emph{a priori} know that the eigenfunction~$\psi$
belongs to $W^{2,2}(\Real^d)$.
This is just because this Sobolev space 
coincides with the domain of~$H_0$ due to the elliptic regularity
(initially, when defining the operator via its sesquilinear form,
we only know that the domain of~$H_0$ consists of functions
$\psi \in W^{1,2}(\Real^d)$ such that $\Delta\psi \in L^2(\Real^d)$).  
This subtlety will become crucial when we deal 
with electromagnetic perturbations below,
allowing critical singularities.
Then an extra regularisation of the multiplier $iT_0\psi$
consists in replacing the gradient in~\eqref{regularised}  
by difference quotients,
as originally proposed in our work~\cite{CK2}.
Altogether, proceeding rigorously with the regularised multiplier
and taking the limits in the right order is rather painful.
This is probably the reason why necessary regularisation schemes
are usually omitted in the literature.

Finally, let us observe that
there is yet another support for the choice~\eqref{dilation},
at least if we deal with the Laplacian and its perturbations.
In fact, the conjugate operator~$T_0$ by itself arises as 
a commutator with the Laplacian:
$$
  T_0 = i\left[H_0, \mbox{$\frac{1}{4}$} x^2\right] .
$$
Consequently,
\begin{equation*}
  \left\langle \mbox{$\frac{1}{4}$} x^2 \right\rangle''
  = \langle T_0 \rangle'
  = \big\langle i[H_0,T_0] \big\rangle
  \,,
\end{equation*}
so the positivity of the commutator $i[H_0,T_0]$
actually shows that the expectation value 
of the square of the magnitude of the position
is a convex function in time: 
there is a \emph{dispersion}. 

\section{Electromagnetic perturbations}\label{Sec.electro}
Of course, the absence of eigenvalues of the free Hamiltonian~$H_0$
can be proved more straight\-for\-ward\-ly 
(\eg, by using the Fourier transform).
However, the advantage of the present method based 
on the virial theorem is that it is much more robust.
In particular, the same conjugate operator~$T_0$ 
applies to electric perturbations of~$H_0$
and its magnetic version enables one to deal 
with magnetic perturbations of~$H_0$, too.

Given a scalar function (electric potential) $V:\Real^d \to \Real$
and a vector-valued function (magnetic potential) $A:\Real^d \to \Real^d$,
consider the electromagnetic Hamiltonian
$$
  H_{A,V} := (-i\nabla -A)^2 + V 
  .
$$
Assume the minimal hypotheses $V \in L_\mathrm{loc}^1(\Real^d)$
and $A \in L_\mathrm{loc}^2(\Real^d)$ 
to give a meaning of the action of $H_{A,V}$ in the sense of distributions.
Moreover, assume that~$|V|$ is relatively form-bounded with respect to 
the magnetic Laplacian $-\Delta_A := (-i\nabla -A)^2$
with the relative bound less than one.
Then $H_{A,V}$ is customarily realised as a self-adjoint operator in $L^2(\Real^d)$
with the form domain of $H_{A,V}$ being the magnetic Sobolev space 
$
  W_A^{1,2}(\Real^d) 
  := \{\psi \in L^2(\Real^d) : \nabla_{\!A}\psi \in L^2(\Real^d)\}
$,
where $\nabla_{\!A} := \nabla -iA$ is the magnetic gradient.
Of course, $H_{0,0} = H_0$ is the free Hamiltonian.

\subsection{Electric perturbations}
In the magnetic-free case, one has 
$$
  i[H_{0,V},T_0] = 2H_0 - x\cdot \nabla V \,,
$$
so the virial identity 
(to be compared with the classical formula~\eqref{virial})
reads
\begin{equation}\label{virial.electric}
  2 \, \|\nabla\psi\|^2  - \int_{\Real^d} x\cdot \nabla V \, |\psi|^2 = 0  
\end{equation}
whenever~$\psi$ is an eigenfunction of $H_{0,V}$.

Clearly, the pointwise repulsivity condition  
\begin{equation}\label{repulsive}
  x\cdot \nabla V \leq 0
\end{equation}
implies a contradiction,
therefore the absence of eigenvalues of $H_{0,V}$.
Less restrictively, it is enough to assume the smallness
of the positive part $(x\cdot \nabla V)_+$ in 
the following integral sense:
\begin{equation}\label{repulsive.bis}
  \exists b < 2 \,, \quad
  \forall \psi \in W^{1,2}(\Real^d)
  \,, \qquad
  \int_{\Real^d} (x\cdot \nabla V)_+ \, |\psi|^2 
  \leq b \int_{\Real^d} \, |\nabla\psi|^2 
  \,.
\end{equation}
In order to justify~\eqref{virial.electric} via the method of multipliers, 
our regularisation scheme described in Section~\ref{Sec.free}
requires the extra regularity condition
\begin{equation}\label{extra}
  V \in W_\mathrm{loc}^{1,p}(\Real^d) 
  \,,
  \qquad \mbox{where} \qquad
  p 
  \begin{cases}
    = 1 & \mbox{if} \quad d=1 \,, \\
    > 1 & \mbox{if} \quad d=2 \,, \\
    d/2 & \mbox{if} \quad d \geq 3 \,. \\
  \end{cases}
\end{equation}

The repulsivity condition~\eqref{repulsive.bis} 
can be replaced by the following smallness condition,
in which case~\eqref{extra} is not needed:
\begin{equation}\label{small}
 \exists b < \frac{2}{d+2} \,, \quad
  \forall \psi \in W^{1,2}(\Real^d)
  \,, \qquad
\begin{aligned}
  \int_{\Real^d} |V| \, |\psi|^2 
  \leq b \int_{\Real^d} \, |\nabla\psi|^2 \,,
  \\
  \int_{\Real^d} |x|^2 \, |V|^2 \, |\psi|^2 
  \leq b^2 \int_{\Real^d} \, |\nabla\psi|^2 \,.
\end{aligned}  
\end{equation}
Indeed, it is enough to integrate by parts in
the second term on the left-hand side of~\eqref{virial.electric}
and use the Schwarz inequality.
Here~\eqref{extra} can be relaxed because 
the identity obtained after the integration by parts
is actually the initial formula to which one arrives by the method of multipliers, 
so differentiating~$V$ is not needed. 

Let us summarise the obtained results into the following theorem.
\begin{theorem}
Assume~\eqref{repulsive.bis} or~\eqref{small}.
In the former case assume in addition~\eqref{extra}.
Then $\sigma_\mathrm{p}(H_{0,V}) = \varnothing$.
\end{theorem}

This theorem is a very special case of a series of recent results 
obtained in \cite[Thm.~3]{FKV2} and \cite[Thm.~3.4]{CFK}.
However, a first rigorous proof of~\eqref{virial.electric} 
(under alternative regularity hypotheses about~$V$) 
goes back to Weidmann~\cite{Weidmann_1967}.

\subsection{Magnetic perturbations}
When there is a magnetic field,
the conjugate operator~\eqref{dilation}
should be replaced by its magnetic version
\begin{equation*}%\label{dilation.magnetic}
  T_A
  := \frac{x \cdot p_A + p_A \cdot x}{2}
  = -i \, x \cdot \nabla_{\!A} - i \, \frac{d}{2}
  \,,
\end{equation*}
where $p_A := -i\nabla_{\!A}$ is the magnetic momentum.
For simplicity, let us consider purely magnetic 
perturbations of the free Hamiltonian.
Then
$$
  i[H_{A,0},T_A] = 2 H_{A,0} 
  + (x \cdot B) \cdot p_A
  + p_A \cdot (x \cdot B) 
  \,,
$$
where $B := \nabla A - (\nabla A)^T$ is the magnetic tensor.
Consequently, the virial identity reads
$$
  2 \, \|\nabla_{\!A}\psi\|^2 
  + 2 \Im \int_{\Real^d} (x\cdot B)\cdot \psi \overline{\nabla_{\!A}\psi}
  = 0 
$$
whenever~$\psi$ is an eigenfunction of~$H_{A,0}$.

Using the Schwarz inequality, we get a contradiction,
and therefore the absence of eigenvalues of~$H_{A,0}$,
provided that the following smallness condition holds:
\begin{equation}\label{small.magnetic}
  \exists b < 1 \,, \quad
  \forall \psi \in W_A^{1,2}(\Real^d)
  \,, \qquad
  \int_{\Real^d} |x|^2 \, |B|^2 \, |\psi|^2 
  \leq b^2 \int_{\Real^d} \, |\nabla_{\!A}\psi|^2 
  \,.
\end{equation}
Our regularisation scheme 
described in Section~\ref{Sec.free}
requires the extra regularity condition
\begin{equation}\label{extra.magnetic}
  A \in W_\mathrm{loc}^{1,2p}(\Real^d) 
  \,,
\end{equation}
where~$p$ is as in~\eqref{extra}.

We have therefore established the following theorem.
\begin{theorem}\label{Thm.magnetic}
Assume~\eqref{small.magnetic} and~\eqref{extra.magnetic}.
Then $\sigma_\mathrm{p}(H_{A,0}) = \varnothing$.
\end{theorem}

It is physically important that the fundamental hypothesis~\eqref{small.magnetic} 
is gauge invariant (\ie, it does not depend on the choice of~$A$
for a given magnetic field~$B$).

Theorem~\ref{Thm.magnetic} 
is a very special case of a series of recent results 
obtained in \cite[Thm.~3]{FKV2} and \cite[Thm.~3.4]{CFK}.
The sufficient conditions which guarantee the absence 
of eigenvalues of~$H_{A,V}$
follow from a full electromagnetic virial identity there.
A first rigorous implementation of the virial theorem
for magnetic Schr\"odinger operators 
goes back to Kalf~\cite{Kalf_1977}.

\subsection{Low versus high dimensions}
It is interesting that spectral conclusions can be obtained
on the basis of functional inequalities of the type~\eqref{repulsive.bis},
\eqref{small} and~\eqref{small.magnetic}.  
Explicit sufficient conditions 
to verify the functional inequalities in high dimensions $d \geq 3$
follow by the Hardy inequality
\begin{equation}\label{Hardy}
  \forall \psi \in W^{1,2}(\Real^d)
  \,, \qquad
  \int_{\Real^d} |\nabla\psi|^2 
  \geq \left( \frac{d-2}{2} \right)^2
  \int_{\Real^d} \frac{|\psi|^2}{|x|^2} 
  \,.
\end{equation}

On the other hand, 
the smallness condition~\eqref{small} cannot be satisfied for a non-trivial~$V$
in low dimensions $d=1,2$.
This is because of the criticality of the Laplacian in these low dimensions,
meaning that any inequality of the type
$H_0 \geq |V|$ necessarily implies $V=0$.

However, hypothesis~\eqref{small.magnetic} 
(and other sufficient conditions stated in terms of
the magnetic Laplacian~$-\Delta_{A}$)
is non-void even in dimension $d=2$
due to the existence of magnetic Hardy inequalities 
\cite{Laptev-Weidl_1999,CK}.

\section{Non-self-adjoint perturbations}\label{Sec.NSA}
There are recent motivations to consider \emph{complex}
electromagnetic fields, including quantum mechanics
\cite{KSTV,K13}.
It is clear already from the manipulations in~\eqref{contradiction}
that the idea based on the virial theorem becomes useless in this case.
On the other hand, the method of multipliers turns out to be more flexible.

\begin{figure}[h!]
\begin{center}
\includegraphics[width=0.8\textwidth]{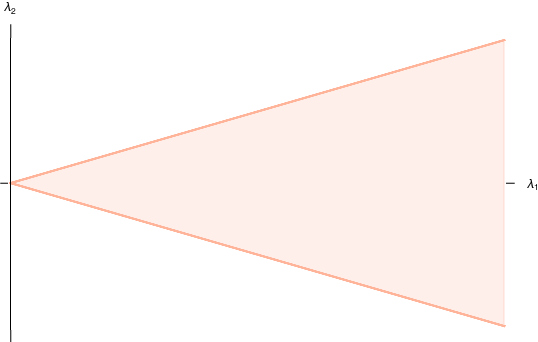}
\end{center}
\caption{The troublesome cone $|\lambda_2| < \lambda_1$.}
\label{Fig.cone}
\end{figure}

Let us demonstrate it on the eigenvalue problem 
for the magnetic-free Hamiltonian
\begin{equation}\label{Helmholtz} 
  H_{0,V}\psi = \lambda\psi \,, 
\end{equation}
where both the potential~$V$ and the eigenvalue~$\lambda$
are allowed to be complex now.
We set $\lambda_1 := \Re\lambda$ and $\lambda_2 := \Im\lambda$,
and analogously for~$V$.
For simplicity, let us assume the following subordination condition:
\begin{equation}\label{subordinate}
  \exists b < 1 \,, \quad
  \forall \psi \in W^{1,2}(\Real^d)
  \,, \qquad
  \int_{\Real^d} (|V_1| + |V_2|) \, |\psi|^2 
  \leq b \int_{\Real^d} \, |\nabla\psi|^2 
  \,.
\end{equation}
Then the numerical range of~$H_{0,V}$
is contained in the cone $|\lambda_2| \leq \lambda_1$,
so it is enough to explore the presence of eigenvalues there,
see Figure~\ref{Fig.cone}.
(In fact, since~$b$ is assumed to be strictly less than one,
the eigenvalues may lie in the interior of the cone only.)

As in Section~\ref{Sec.multiplier}
(and disregarding the necessary regularisation procedures),
take an inner product of both sides of~\eqref{Helmholtz} 
with the function~$iT_0\psi$, where~$T_0$ is given by~\eqref{dilation},
and take twice the real part of the obtained identity.
This leads to the virial identity
\begin{equation}\label{4.36} 
  2 \, \|\nabla\psi\|^2 
  -\int_{\Real^d} x\cdot\nabla V_1 \, |\psi|^2
  -2 \Im(x\cdot\nabla\psi,V_2 \psi)
  = - 2 \, \lambda_2 \, \Im(x\cdot\nabla\psi,\psi) 
  \,,
\end{equation}
which is a non-self-adjoint counterpart of~\eqref{virial.electric}.
Because of the right-hand side of~\eqref{4.36} with no obvious sign, 
this identity by itself does not enable one to conclude 
with the total absence of eigenvalues with $\lambda_2 \not= 0$.
(Interestingly, real eigenvalues can be excluded by conditions
of the type~\eqref{repulsive.bis} and~\eqref{small}.)

The idea of \cite{Barcelo-Vega-Zubeldia_2013}
(which seems to go back to~\cite{Eidus_1962} and~\cite{Ikebe-Saito_1972},
while it was developed in the present spectral setting in~\cite{FKV})
is to compensate the appearance of the right-hand side of~\eqref{4.36} 
for non-real eigenvalues by further identities obtained by 
using different multipliers.
First, taking an inner product of both sides of~\eqref{Helmholtz} 
with the function $\psi$ 
and taking the real part of the obtained identity,
we get
\begin{equation}\label{4.32} 
  \|\nabla\psi\|^2
  + \int_{\Real^d} V_1 \, |\psi|^2
  = \lambda_1 \,  \|\psi\|^2
  \,.
\end{equation}
Second, taking an inner product of both sides of~\eqref{Helmholtz} 
with the function $|x|\psi$ 
and taking the real and imaginary part of the obtained identity,
we respectively get
\begin{equation}\label{4.33} 
  \int_{\Real^d} |x| \, |\nabla\psi|^2
  -\frac{d-1}{2} \int_{\Real^d} \frac{|\psi|^2}{|x|}
  + \int_{\Real^d} |x| \, V_1 \, |\psi|^2
  = \lambda_1 \,  \int_{\Real^d} |x| \, |\psi|^2
\end{equation}
and 
\begin{equation}\label{4.35} 
  \Im \int_{\Real^d} \frac{x}{|x|} \cdot \overline{\psi} \nabla\psi 
  + \int_{\Real^d} |x| \, V_2 \, |\psi|^2
  = \lambda_2 \,  \int_{\Real^d} |x| \, |\psi|^2
  \,.
\end{equation}
By taking the clever sum
$$
  \eqref{4.36} - \eqref{4.32} 
  + \frac{|\lambda_2|}{\sqrt{\lambda_1}} \, \eqref{4.33} 
  - 2\, \sqrt{\lambda_1} \, \sgn(\lambda_2) \, \eqref{4.35} 
  \,,
$$  
we arrive at the ultimate identity
\begin{align}\label{ultimate} 
\lefteqn{
  \|\nabla\psi^-\|^2
  + \frac{|\lambda_2|}{\sqrt{\lambda_1}}
  \int_{\Real^d} |x| \left( 
  |\nabla\psi^-|^2 - \frac{d-1}{2} \, \frac{|\psi|^2}{|x|}
  \right)
  }  
  \\
  &
  + (d-1) \int_{\Real^d} V_1 \, |\psi|^2
  + \frac{|\lambda_2|}{\sqrt{\lambda_1}} \int_{\Real^d} |x| \, V_1 \, |\psi|^2
  + 2 \, \Re \int_{\Real^d} V \, x \cdot \psi^- \overline{\nabla\psi^-}
  = 0
  \nonumber
  \,,
\end{align}
where 
$$
  \psi^-(x) := \E^{-i \sqrt{\lambda_1} \, \sgn(\lambda_2) \, |x|}
  \, \psi(x)
  \,.
$$
No condition of the type~\eqref{extra} is needed to justify this result.

Various sufficient conditions for the absence of eigenvalues of~$H_{0,V}$
can be derived from~\eqref{ultimate}. 
This has been done in a series of recent papers \cite{FKV,FKV2,CFK},
including the magnetic field. 

For instance, let $d\geq 3$, so that the second term 
on the first line of~\eqref{ultimate} is non-negative 
by a weighted Hardy inequality, 
and assume for simplicity that the potential~$V$ is purely imaginary. 
Then~$H_{0,V}$ has no eigenvalues in the cone $|\lambda_2| \leq \lambda_1$
provided that the following condition holds:
\begin{equation}\label{small.imaginary}
  \exists b < 1/2 \,, \quad
  \forall \psi \in W^{1,2}(\Real^d)
  \,, \qquad
  \int_{\Real^d} |x|^2 \, |\Im V|^2 \, |\psi|^2 
  \leq b^2 \int_{\Real^d} \, |\nabla\psi|^2 
  \,.
\end{equation}
Let us summarise the sufficient conditions which guarantee 
the total absence of eigenvalues of~$H_{0,V}$ in this special case
into the following theorem.

\begin{theorem} 
Let $d \geq 3$ and $\Re V = 0$.
Assume conditions~\eqref{small.imaginary} and~\eqref{subordinate}. 
Then $\sigma_\mathrm{p}(H_{0,V}) = \varnothing$.
\end{theorem}

\section{Uniform resolvent estimates}\label{Sec.uniform}
The power of the method of multipliers 
as developed in the preceding section 
can be used to derive finer spectral properties, 
going beyond the mere absence of eigenvalues.
As an example, let us use it to derive 
uniform resolvent estimates for the free Hamiltonian~$H_0$.

Because of the fundamental identity
$$
  \forall \lambda \not\in \sigma(H_0) = [0,+\infty)
  \,, \qquad
  \| (H_0-\lambda)^{-1} \| 
  = \frac{1}{\dist(\lambda,\sigma(H_0))}
  \,,
$$
there is no hope to bound the operator norm 
of the resolvent by a constant independent of~$\lambda$.
However, it is well known that it is possible 
when the resolvent is reconsidered as an operator
acting between different spaces.
An example of such \emph{uniform resolvent estimates}
is the celebrated result of Kato and Yajima
\cite[Thm.~1]{Kato-Yajima_1989}:
\begin{Theorem}\label{Thm.uniform}
Let $d \geq 3$. Then
$
\displaystyle
  \sup_{\lambda \not\in \sigma(H_0)} 
  \big\| |x|^{-1} (H_0-\lambda)^{-1} |x| \big\| 
  =: C < \infty
$.
\end{Theorem}
This result implies not only that~$H_0$ possesses no eigenvalues
but that its spectrum is actually purely absolutely continuous. 

\begin{proof}[Theorem~\ref{Thm.uniform} via the method of multipliers]
Consider the resolvent equation 
\begin{equation}\label{resolvent}
  (H_0-\lambda)\psi = f 
  \,
\end{equation}
where $\lambda \not\in \sigma(H_0)$ and $f \in C_0^\infty(\Real^d)$. 

Taking an inner product of both sides of~\eqref{resolvent}
with~$\psi$ and taking the real part of the obtained identity,
we get 
\begin{equation}\label{id1}
  \|\nabla\psi\|^2 - \lambda_1 \, \|\psi\|^2 
  = \Re (\psi,f) 
  \,.
\end{equation}
Consequently, if $\lambda_1 \leq 0$, 
then the Hardy inequality~\eqref{Hardy}
and the Schwarz inequality imply
$
  \big\||x|^{-1} \psi\big\|  
  \leq  
  C_1 \,
  \big\||x|f\big\|
$
with
$$
  C_1 := \left( \frac{2}{d-2} \right)^2
  \,.
$$
This is the desired resolvent estimate 
in the complex half-plane $\lambda_1 \leq 0$.

Similarly, taking an inner product of both sides of~\eqref{resolvent}
with~$\psi$ and taking the imaginary part of the obtained identity,
we get 
\begin{equation}\label{id2}
  - \lambda_2 \, \|\psi\|^2 
  = \Im (\psi,f) 
  \,.
\end{equation}
Summing up~\eqref{id1} with the absolute value of~\eqref{id2} 
and using the Hardy and Schwarz inequalities as above, 
we get 
$
  \big\||x|^{-1} \psi\big\|
  \leq  
  C_2 \,
  \big\||x|f\big\|
$
with
$$
  C_2 := \sqrt{2} \, \left( \frac{2}{d-2} \right)^2
$$
in the complex region $0 < \lambda_1 \leq |\lambda_2|$.

It remains to analyse the troublesome cone $|\lambda_2| < \lambda_1$,
see Figure~\ref{Fig.cone}.	
Proceeding as in Section~\ref{Sec.NSA}, we arrive at the identity
\begin{align}\label{ultimate.bis} 
\lefteqn{
  \|\nabla\psi^-\|^2
  + \frac{|\lambda_2|}{\sqrt{\lambda_1}}
  \int_{\Real^d} |x| \left( 
  |\nabla\psi^-|^2 - \frac{d-1}{2} \, \frac{|\psi|^2}{|x|}
  \right)
  }  
  \\
  &
  = (d-1) \Re \int_{\Real^d} \overline{\psi} \, f
  + \frac{|\lambda_2|}{\sqrt{\lambda_1}} \int_{\Real^d} |x| \, \overline{\psi} \, f
  + 2 \, \Re \int_{\Real^d} f^- \, x \cdot \overline{\nabla\psi^-}
  \nonumber
  \,.
\end{align}
Indeed, this identity coincides with~\eqref{ultimate}
after the formal identification $f = -Vu$. 
The first term on the right-hand side of~\eqref{ultimate.bis} 
is estimated by means of the Hardy inequality~\eqref{Hardy} as follows:
$$
\begin{aligned}
  \left| (d-1) \Re \int_{\Real^d} \overline{\psi} \, f \right|
  &\leq  (d-1) \, \big\||x|^{-1} \psi^-\big\| \, \big\||x|f\big\|
  \leq 2 \, \frac{d-1}{d-2} \, \big\|\nabla \psi^-\big\| \, \big\||x|f\big\|
  \\
  &\leq \delta \, \big\|\nabla \psi^-\big\|^2 
  + \frac{1}{\delta} \, \left(\frac{d-1}{d-2}\right)^2 \, \big\||x|f\big\|^2
  \,,
\end{aligned}
$$
where the last inequality holds with any $\delta \in (0,1)$.
The last term on the right-hand side of~\eqref{ultimate.bis} 
is also easy to estimate:
$$
\begin{aligned}
  \left| 2 \, \Re \int_{\Real^d} f^- \, x \cdot \overline{\nabla\psi^-}  \right|
  &\leq 2 \, \big\|\nabla \psi^-\big\| \, \big\||x|f\big\|
  \leq  \delta \, \big\|\nabla \psi^-\big\|^2 
  + \frac{1}{\delta} \, \big\||x|f\big\|^2
  \,.
\end{aligned}
$$
To estimate the middle term on the right-hand side of~\eqref{ultimate.bis},
we note that~\eqref{id2} implies
$$
  \|\psi^-\|^2 \leq \frac{1}{|\lambda_2|} \, 
  \big\||x|^{-1} \psi^-\big\| \, \big\||x|f\big\|
  \,.
$$
Consequently, using that we are inside the cone $|\lambda_2| < \lambda_1$,
$$
\begin{aligned}
  \left| \frac{|\lambda_2|}{\sqrt{\lambda_1}} 
  \int_{\Real^d} |x| \, \overline{\psi} \, f \right|
  &\leq  \frac{|\lambda_2|}{\sqrt{\lambda_1}} \,
  \big\|\psi^-\big\| \, \big\||x|f\big\|
  \leq \big\||x|^{-1} \psi^-\big\|^{1/2} \, \big\||x|f\big\|^{3/2}
  \\
  &\leq \sqrt{\frac{2}{d-2}} \,
  \big\|\nabla\psi^-\big\|^{1/2} \, \big\||x|f\big\|^{3/2}
  \\
  &\leq \frac{1}{4} \, \delta \, \big\|\nabla\psi^-\big\|^2
  + \frac{3}{4} \, \frac{1}{\delta^{1/3}} \,
  \left(\frac{2}{d-2}\right)^{2/3}
  \, \big\||x|f\big\|^{2}
  \,.
\end{aligned}
$$
Neglecting the second term on the left-hand side of~\eqref{ultimate.bis}
(which is non-negative whenever $d \geq 3$ by a weighted Hardy inequqlity),
we thus arrive at the inequality
$$
  \|\nabla\psi^-\|^2 \, 
  \left(1 - \delta - \mbox{$\frac{1}{4}$} \delta -  \delta \right)
  \leq \big\||x|f\big\|^2 \,
  \left[
  \frac{1}{\delta}
  \left(\frac{d-1}{d-2}\right)^2
  + \frac{3}{4} \, \frac{1}{\delta^{1/3}} \,
  \left(\frac{2}{d-2}\right)^{2/3}
  + \frac{1}{\delta}
  \right] .
$$
Estimating $\delta^{-1/3} > \delta^{-1}$,
optimising with respect to~$\delta$
(\ie, choosing $\delta:= 2/9$)
and using the Hardy inequality~\eqref{Hardy} once more,
we eventually get 
$
\big\||x|^{-1} \psi\big\|
  \leq  
  C_3 \,
  \big\||x|f\big\|
$
with
$$
  C_3 := \frac{6}{d-2} \left[
  \left(\frac{d-1}{d-2}\right)^2
  + \frac{3}{4} \,  
  \left(\frac{2}{d-2}\right)^{2/3}
  + 1
  \right]^{1/2} 
$$
in the cone $|\lambda_2| < \lambda_1$.

In summary, since $C_1 < C_2 < C_3$, 
we have got
$
\big\||x|^{-1} \psi\big\|
  \leq  
  C_3 \,
  \big\||x|f\big\|
$  
with $C \leq C_3$ in the whole resolvent set
$\lambda \in \Com \setminus [0,\infty)$.	
\hfill\qed 
\end{proof}

\section{Relativistic operators}\label{Sec.Dirac}
The approach described in Section~\ref{Sec.NSA}
can be adapted to electromagnetic Schr\"o\-ding\-er 
operators with \emph{matrix-valued} potentials.
This has been done in~\cite{CFK}, 
where we also applied the results to establish the absence
of eigenvalues of Pauli and Dirac operators. 
While we consider arbitrary dimensions in~\cite{CFK}, 
let us focus on dimension $d=3$ to present our results in a succinct way.
For simplicity, let us also restrict to the purely magnetic case, \ie\ $V=0$.

Given a locally square integrable potential $A:\Real^3 \to \Real^3$ 
as in Section~\ref{Sec.electro},
consider the (self-adjoint) Dirac operator
$$
  D_{A,0} := - i \alpha \cdot \nabla_{\!A} + m \, \alpha_0
$$ 
in the Hilbert space $L^2(\Real^3)^4$ with domain $W_A^{1,2}(\Real^3)^4$.
Here $m \geq 0$ is the mass of the particle 
and $\alpha_0, \alpha_1, \alpha_2, \alpha_3$ are the standard $4 \times 4$ 
Hermitian Dirac matrices satisfying the anticommutation rules 
$
  \alpha_j \alpha_k + \alpha_k \alpha_j = 2 \delta_{jk} I_{\Com^4}  
$
with $j,k \in \{0,1,2,3\}$.

The virial theorem and the method of multipliers 
do not seem to apply directly to the Dirac operators,
because of the lack of positivity of certain commutators. 
(However, partial results can be derived \cite[Thm.~4.2.1]{Thaller}.)
Our strategy is to employ the well-known supersymmetric structure instead,
meaning that 
\begin{equation}\label{super}
  {D_{A,0}}^2 = 
  \begin{pmatrix}
    P_{A,0} + m^2 I_{\Com^2}  & 0 \\
    0 &  P_{A,0} + m^2 I_{\Com^2} 
  \end{pmatrix}
\end{equation}
where $P_{A,0} := H_{A,0} \, I_{\Com^2}  + \sigma \cdot B$ 
with $B := \curl A$
is the (self-adjoint) Pauli operator in $L^2(\Real^3)^2$ 
with domain $W_A^{1,2}(\Real^3)^2$.
Here $\sigma_1, \sigma_2, \sigma_3$ are the standard 
$2 \times 2$ Hermitian Pauli matrices
satisfying the anticommutation rules 
$
  \sigma_j \sigma_k + \sigma_k \sigma_j = 2 \delta_{jk} I_{\Com^2} 
$
with $j,k \in \{1,2,3\}$.

As a consequence of~\eqref{super}, the spectral relationship
$\sigma(D_{A,0})^2 = \sigma(P_{A,0}) + m^2$ holds true.
What is more, $\sigma_\mathrm{p}(D_{A,0}) = \varnothing$
if, and only if, $\sigma_\mathrm{p}(P_{A,0}) = \varnothing$.
The proof of the absence of eigenvalues of $D_{A,0}$ 
thus reduces to the absence of eigenvalues of
electromagnetic Schr\"odinger operators $H_{A,V}$,
like in Section~\ref{Sec.electro} but with the generalisation that
the electric potential is allowed to be matrix-valued 
(in fact, $V = \sigma \cdot B$ in the Pauli case).
In~\cite{CFK}, we also allow for possibly non-Hermitian potentials. 

As one example of the results established in~\cite{CFK}, 
let us present the following theorem (\cf~\cite[Thm.~1.3]{CFK}).  

\begin{theorem} 
Let $d = 3$ and $A \in W_\mathrm{loc}^{1,3}(\Real^3)$.
Assume condition~\eqref{small.magnetic} with $b < 1/14$. 
Then $\sigma_\mathrm{p}(D_{A,0}) = \varnothing$.
\end{theorem}

\section{Boundary perturbations}\label{Sec.boundary}
The flexibility of the method of multipliers,
particularly due to the developments described in Section~\ref{Sec.NSA},
enables one to consider elliptic operators constrained 
to \emph{subdomains} of the Euclidean space.

In~\cite{CK2}, we developed the method to study spectral properties
of the Laplacian $-\Delta_\alpha$
in the half-space $\Omega := \Real^{d-1} \times (0,\infty)$,
subject to Robin boundary conditions 
$$
  -\frac{\partial\psi}{\partial x_d} + \alpha \, \psi = 0 
  \qquad \mbox{on} \qquad
  \partial\Omega = \Real^{d-1} \times \{0\}
  \,,
$$
where $\alpha : \partial\Omega \to \Com$
plays the role of a strongly localised potential.
Under the hypothesis $\alpha \in L^\infty(\partial\Omega)$,
the Robin Laplacian $-\Delta_\alpha$ 
can be realised as an m-sectorial operator in $L^2(\Omega)$,
with form domain $W^{1,2}(\Omega)$.

As one example of the results established in~\cite{CK2}, 
let us present the following theorem (\cf~\cite[Thm.~1.3]{CFK}). 
\begin{theorem} 
Let $\alpha \in W_\mathrm{loc}^{1,\infty}(\partial\Omega)$ be real-valued.
Assume 
\begin{equation}\label{Robin}
  \alpha \geq 0
  \qquad \mbox{and} \qquad
  x \cdot \nabla \alpha \leq 0
  \,.
\end{equation}
Then $\sigma_\mathrm{p}(-\Delta_\alpha) = \varnothing$.
\end{theorem}

Of course, \eqref{Robin}~plays the role of the repulsivity 
condition~\eqref{repulsive}. 
In addition to complex-valued boundary conditions,
we also derive uniform resolvent estimates in~\cite{CK2}.

The half-space can be regarded as a degenerate situation 
of conical domains intensively studied in recent years.
In this respect, let us particularly mention
the proof of the absence of eigenvalues of the Laplacian
in non-convex conical sectors,
subject to no specific boundary conditions
\cite{Bonnet-Fliss-Hazard-Tonnoir_2016}.
On the other hand, it is easy to construct square-integrable solutions 
to the eigenvalue problem in a half-space.

A variant of the method of multipliers for the Dirichlet Laplacian
in repulsive waveguides has been developed in~\cite{D'Ancona-Racke_2012}.
Finally, as another flexibility of the method of multipliers,
let us mention its recent developments 
for the Lam\'e operator in elasticity~\cite{Cossetti_2017}
and polyharmonic operators~\cite{CFK2}.

\subsection*{Acknowledgments}
The first author (L.C.) was supported 
by the grant Ram\'on y Cajal RYC2021-032803-I funded by
MCIN/AEI/10.13039/50110 0011033 and by the European Union
NextGenerationEU/PRTR and by Ikerbasque.
The second author (D.K.) was supported
by the EXPRO grant No.~20-17749X 
of the Czech Science Foundation.

%--------------%
% BIBLIOGRAPHY %
%--------------%
%
%\newpage
%\addcontentsline{toc}{section}{References}
%\bibliography{bib}
%\bibliographystyle{amsplain}
%
 
\providecommand{\bysame}{\leavevmode\hbox to3em{\hrulefill}\thinspace}
\providecommand{\MR}{\relax\ifhmode\unskip\space\fi MR }
\providecommand{\MRhref}[2]{%
  \href{http://www.ams.org/mathscinet-getitem?mr=#1}{#2}
}
\providecommand{\href}[2]{#2}

\end{document}